\documentclass[11pt,reqno,twoside]{amsart}
\usepackage{amssymb,amsmath,amsthm,soul,color}
\usepackage{t1enc}
\usepackage[cp1250]{inputenc}
\usepackage{a4,indentfirst,latexsym}
\usepackage{graphics}
\usepackage{mathrsfs}
\usepackage{cite,enumitem,graphicx}

\usepackage[colorlinks=true,urlcolor=blue,citecolor=red,linkcolor=blue,linktocpage,pdfpagelabels,
bookmarksnumbered,bookmarksopen]{hyperref}

\usepackage[english]{babel}
\usepackage[left=2.61cm,right=2.61cm,top=2.72cm,bottom=2.72cm]{geometry}

\usepackage[metapost]{mfpic}

\usepackage[colorinlistoftodos]{todonotes}
\usepackage[normalem]{ulem}

\allowdisplaybreaks

\makeatletter
\providecommand\@dotsep{5}
\def\listtodoname{List of Todos}
\def\listoftodos{\@starttoc{tdo}\listtodoname}
\@namedef{subjclassname@2020}{\textup{2020} Mathematics Subject Classification}
\makeatother

\newcounter{stepnum}

\newtheorem{thm}{Theorem}[section]
\newtheorem{pro}[thm]{Proposition}
\newtheorem{lem}[thm]{Lemma}

\newcommand{\RT}{{\mathbb{R}^3}}

\newcommand{\D}{D^{1,q}(\mathbb{R}^3)}

\newcommand{\INT}{\int_{\mathbb{R}^3}}

\numberwithin{equation}{section}

\begin{document}

\title[ \MakeLowercase{($p,q$)}-Schr\"odinger-Poisson  system]{Infinitely many solutions for \MakeLowercase{($p,q$)}-Schr\"odinger-Poisson system with concave and convex nonlinearities}
	\author[Y. Du]{Yao Du}
	\author[J. Peng]{Jiahao Peng}
	\address[Y. Du]{\newline\indent
		School of  Science,
		\newline\indent
		XiHua University,
		\newline\indent
		Chengdu 610039, People's Republic of China}
	\email{\href{1052591976@qq.com}{1052591976@qq.com}}
	
	\address[J. Peng]{\newline\indent
		School of  Science,
		\newline\indent
		XiHua University,
		\newline\indent
		Chengdu 610039, People's Republic of China}
	\email{\href{1503379504@qq.com}{1503379504@qq.com}}

	\date{\today}

	\maketitle

 \begin{abstract}
In this paper, we obtain infinitely many solutions for  a  class of quasilinear Schr\"{o}dinger-Poisson system  which is coupled by a Schr\"{o}dinger equation of $p$-Laplacian and a Poisson equation of $q$-Laplacian, involving with concave and convex nonlinearities and indefinite weighted functions.
 \end{abstract}

\section{Introduction}
In this paper, we consider the following Schr$\ddot{\mathrm{o}}$dinger-Poisson system
\begin{eqnarray}\label{equ}
\left \{\begin{array}{ll}
\displaystyle -\Delta_p u-\lambda g(x)|u|^{p-2}u+a(x)\phi |u|^{m-2}u=k(x)|u|^{r-2}u-h(x)|u|^{s-2}u,&\mathrm{in} \ \mathbb{R}^3,\\[1em]
\displaystyle -\Delta_q \phi =a(x) |u|^m,  &\mathrm{in}\ \mathbb{R}^3,\\ [1em]
u\in E,\phi\in D^{1,q}(\mathbb{R}^3),
\end{array}
\right.
\end{eqnarray}
where $\Delta_i u=\hbox{div}(|\nabla u|^{i-2}\nabla u)\ (i=p,q)$, and
\begin{eqnarray}
1<r<p<s<+\infty,\ p<3,\   r<\frac{qm}{q-1}\label{n1.2}
\end{eqnarray}
with $q$ and $m$ fulfilling
\begin{eqnarray}\label{2.6}
 \max \left\{1,\frac{3p}{5p-3}\right\}<q<3, \ 1<m<\frac{(q^*-1)p^*}{q^*}.
\end{eqnarray}
The Banach space $E$  is defined as the completion of  $C^{\infty}_0(\RT)$ under the norm
\begin{eqnarray*}
    \|u\|:=\|u\|_{D^{1,p}(\RT)}+|u|_{s,h}=\left(\INT | \nabla  u|^p dx\right)^{1/p}+\left(\INT h(x)|u|^s dx \right)^{1/s}.
\end{eqnarray*}
For $1\leq t<3$ and  $t^*:=\frac{3t}{3-t}$, see \cite[Definition 7.2.1]{2013Willem}, on the function space
$$D^{1,t}(\mathbb{R}^3)=\left\{u\in L^{t^*}(\RT):\ \nabla u\in L^t(\RT, \RT)\right\},$$
 we define the  norm
\begin{eqnarray*}
\|u\|_{D^{1,t}(\RT)}:=\left(\INT|\nabla u|^t dx\right)^{\frac{1}{t}}.
\end{eqnarray*}
The embedding, see \cite[Proposition 7.2.2]{2013Willem},
\begin{eqnarray}
D^{1,t}(\RT)\hookrightarrow L^\ell(\RT)\label{1.3}
\end{eqnarray}
holds if and only if $\ell=t^*$, where $t^*$ is called   the critical Sobolev exponent of the embedding \eqref{1.3} and $L^\ell(\RT)$ denotes the Lebesgue space for $1\le \ell<\infty$, equipped with the norm
\begin{eqnarray*}
    \|u\|_\ell:=\left(\INT |u|^\ell dx\right)^{\frac{1}{\ell}}.
\end{eqnarray*}
 We also define the best Sobolev constant
\begin{eqnarray}\label{38}
S_t:=\inf _{u \in D^{1,t}(\RT) \backslash\{0\}} \frac{\|u\|_{D^{1,t}(\RT)}}{\|u\|_{t^{*}}^{t}}.
\end{eqnarray}

Motivated by the studies in \cite{2023DuSuWang,2003LiuLi}, the problem \eqref{equ} under consideration is a nonlinear elliptic system coupling a Schr\"{o}dinger equation with $p$-Laplacian and a Poisson equation with $q$-Laplacian, including concave and convex nonlinearities. The $p$-Laplacian operator arises in nonlinear fluid dynamics, where the value of $p$ depends on flow velocity and material properties. For further details on its physical origins, see \cite{2008BenediktGrigKotrlaTakac}.

The classical Schr\"{o}dinger-Poisson system (corresponding to $p=q=m=2$ in \eqref{equ}) appears naturally in quantum mechanics and semiconductor theory, see  \cite{1981BenguriaBrezisLieb,1990MarkwichRinghoferSchmeiser}.  It has been extensively studied via variational methods, we refer to  \cite{2008Ambrosetti,2008AR,1998BenciFortunato, 2008AA,2004DM-1,2004DM-2,2002d'Avenia,2014LG,2006Ruiz,2009ZZ,2018ZT}  and references therein. Recently, several works have begun to explore systems involving $p$-Laplacian. For the quasilinear Schr\"odinger-Poisson system   we refer to \cite{2021DuSuWang,2022DuSuWang,2022DuSuWang505,2024DuSuWang, 2023DuSuWang}. When $q=2$, the first author of this paper, together with Su and Wang \cite{2021DuSuWang, 2022DuSuWang, 2022DuSuWang505}, established the existence of nontrivial solutions and ground state solutions for   quasilinear Schr\"odinger-Poisson system.

The case $q \ne 2$ introduces three main obstacles to the system \eqref{equ}. The first one is that the lack of linearity for operator in the Poisson equation prevents the direct application of the Lax-Milgram theorem to obtain a unique solution $\phi_u \in D^{1,q}(\mathbb{R}^3)$ for each $u \in E$.  In \cite{2023DuSuWang}, the first author, Su and Wang used the Minty-Browder theorem to overcome this obstacle, which is essential for constructing the variational framework.   The second one is that the solution $\phi_u$ lacks an explicit expression, making it difficult to derive further properties.  Based on the uniqueness of $\phi_u$, the key properties of the solution $\phi_u$ for the second equation in \eqref{equ} with $a \equiv 1$ have been established in  \cite{2023DuSuWang}. Third, the case $q \ne 2$ complicates the variational structure. For $q=2$, the energy functional associated with quasilinear system is well-defined and of class $C^1$ via a standard reduction procedure (see \cite{2021DuSuWang, 2022DuSuWang, 2022DuSuWang505}). For $q \ne 2$, this is a new challenge. Inspired by \cite{2010Yu}, the authors of \cite{2023DuSuWang} used methods from nonlinear functional analysis to prove G\^{a}teaux differentiability and the continuity of the derivative, showing that the functional is of class $C^1$.   The first author of this paper and Su  in \cite{2024DuSuWang} proved the existence of nontrivial solutions of quasilinear Schr\"odinger-Poisson system with Berestycki-Lions type conditions.

For semilinear elliptic equations with concave and convex nonlinearities, we refer to   \cite{1994AmbrosettiBrezisCerami,1995BartschWillem,2006Wu,2010Wu,2003LiuLi} and  related references.  The classical  Schr\"{o}dinger-Poisson system with  concave and convex nonlinearities has been studied in \cite{2018MaoShao,2021SunWu,2017LiTang,2021ChenLi,2017LeiSuo}. Do the quasilinear Schr\"odinger-Poisson system involving with concave and convex nonlinearities have solutions? In the present paper, we establish the infinitely many solution of the quasilinear  Schr\"odinger-Poisson system \eqref{equ}.

For a given number $r\in(1,p^*)$, we always denote
\begin{eqnarray}\label{456}
r_0:=\frac{3p}{3p-r(3-p)}.
\end{eqnarray}

We impose the following conditions on the weighted functions:\\
(a) $a(x)\in L^{\infty }(\mathbb{R}^3)\cap L^{\gamma_{0} }(\mathbb{R}^3),a(x)\ge 0,$ where $\gamma_{0}=\frac{p^*}{p*-\gamma }$ and $\gamma=\frac{mq^*}{q*-1}$; \\
(g) $ g(x)\in L^{3/p} (\mathbb{R}^3)\cap L^{\infty } (\mathbb{R}^3)$ is indefinite in sign;\\
(k) $ k(x)\in L^{r_{0} }(\mathbb{R}^3) \cap L^{\infty } (\mathbb{R}^3) ,k(x)\ge 0$ and $ k(x)\not\equiv 0$;\\
(h) $ h(x)\in L^{\infty}(\mathbb{R}^3)$, $h(x)\ge 0$\ for a.e. $ x\in \mathbb{R}^3 $.\\

Let $\lambda_{-1}$ and $\lambda_1$  denote the principal negative and positive eigenvalues of the eigenvalue problem (see \cite{1997DrabekHuang,1995AllegrettoHuang})
\begin{eqnarray}
\left \{\begin{array}{ll}
\displaystyle -\Delta_p u=\lambda g(x)|u|^{p-2}u,\\
u\in D^{1,p}( \RT).
\end{array}
\right. \label{101}
\end{eqnarray}

Our main result is as follows.
\begin{thm}\label{thm1.1}
Assume that \eqref{n1.2}, \eqref{2.6}, and \textup{(a), (g), (k), (h)} hold. Then, for every $\lambda \in (-\lambda_{-1}, \lambda_1)$, problem \eqref{equ} admits infinitely many solutions with negative energy.
\end{thm}

Under  assumptions (g), (k), (h) on the weighted functions $g(x)$, $k(x)$, $h(x)$,  Liu and Li in \cite[Theorem  1.1]{2003LiuLi} established the infinitely many solutions for the $p$-Laplacian problem, such problem is $a\equiv0$ in system \eqref{equ}. For the Schr\"{o}dinger-Poisson system with $(p,q)$-Laplacian, several fundamental difficulties arise.  As described in \cite{2023DuSuWang}, these include: the inability to apply the Lax-Milgram theorem directly; the lack of an explicit expression for $\phi_u$; and the need to prove that the variational functional is $C^1$. Additionally, the presence of concave and convex nonlinearities and indefinite weighted functions, compounded by the non-reflexivity of the function space $E$,   significantly impact the existence and multiplicity of solutions. Theorem \ref{thm1.1} extends the results of \cite{2003LiuLi, 2023DuSuWang}.

This paper is organized as follows: in Section \ref{sec2}, we introduce the function spaces and preliminary results. Section \ref{sec3} is devoted to constructing the variational framework and proving the $C^1$-regularity of the energy functional. Finally, in Section \ref{sec4}, we prove Theorem \ref{thm1.1} by verifying coercive  and the Palais-Smale condition  of the energy functional.

\section{Preliminaries}\label{sec2}
The following function spaces and notations will be employed in our analysis of problem \eqref{equ}.
\begin{itemize}
\item $L^r_l(\mathbb{R}^3)$ is the weighted Lebesgue space consisting of all measurable functions $u$ such that
 \begin{eqnarray*}
\INT l(x)|u|^rdx<\infty
 \end{eqnarray*}
with a given nonnegative measurable weighted function $l$ and a real number $r>1$.  Associated with this space is the seminorm
\begin{eqnarray*}
|u|_{r,l}:=\left(\INT l(x)|u|^r dx\right)^{\frac{1}{r}}.
\end{eqnarray*}
\item $L^q_{loc}(\RT)$ denotes the local Lebesgue space for  $1\le q<\infty $.  A function $f$ belongs to this space if for every compact set $K \subset \mathbb{R}^3$,   $\int_{K}|f(x)|^q<\infty$.

    \item $L^\infty(\RT)$ denotes the Banach space of measurable functions that are essentially bounded, equipped with the norm
\begin{eqnarray*}
    \displaystyle\|u\|_\infty:=\mathop{\rm ess\ sup}\limits_{x\in \RT}|u(x)|.
\end{eqnarray*}

\item We denote by $X^{'}$ the dual space of a Banach space $X$, equipped with the norm $\|\cdot\|_{X^{'}}$. The duality pairing between $X^{'}$ and $X$ is denoted by $\langle \cdot, \cdot \rangle$.

\item  Let $C$ be a positive constant that may change from line to line.
\end{itemize}

The reflexivity of a space is an effective tool for obtaining a weakly convergent subsequence. However, we do not know whether $E$ is reflexive.

It is easy to see that the following inclusions hold:
\begin{eqnarray}
    C^{\infty}_0(\RT)\subset E\subset D^{1,p}(\RT)\cap L^s_h(\RT).\label{n2.1}
\end{eqnarray}

The characteristics of the space $E$ are examined by the following lemmas.
\begin{lem}\textup{\cite[Lemma 2.1]{2003LiuLi}}\label{lem2.1}
For $\gamma \in (1,p^*)$, if $l \in L^{\gamma_0}(\mathbb{R}^N) \cap L^{\infty}(\mathbb{R}^N)$ with $l \ge 0$, where $\gamma_0$ is defined in \eqref{456}, then the embedding $D^{1,p}(\mathbb{R}^N) \hookrightarrow L^{\gamma}_l(\mathbb{R}^N)$ is compact. Moreover, the embedding $E \hookrightarrow L^{\gamma}_l(\mathbb{R}^N)$ is also compact.
\end{lem}

\begin{lem}\textup{\cite[Lemma 2.2]{2003LiuLi}}\label{lem2.2} Assume \eqref{n1.2} and (h) hold, then
\begin{description}
  \item[(i)] there exists a constant $C>0$ such that for all $u \in E$,
\begin{eqnarray*}
\frac{1}{p}\INT | \nabla u |^pdx+\frac{1}{s}\INT h(x)|  u |^sdx\le C(\|u\|^p+\|u\|^s).
\end{eqnarray*}
  \item[(ii)] for any $\alpha, \beta > 0$, there exists $c > 0$ such that for every $u \in E$,
 \begin{eqnarray*}
\alpha \INT|\nabla u|^{p}dx+\beta \INT h(x)|u|^{s}x \ge\left\{\begin{array}{ll}
c\|u\|^{s} & \text { if }\|u\| \le 1, \\
c\|u\|^{p} & \text { if }\|u\| \ge 1 .
\end{array}\right.
\end{eqnarray*}
\end{description}
\end{lem}

\section{Variational Framework}\label{sec3}
Be mindful of the conditions $1<q<3$ and $1<m<\frac{(q^*-1)p^*}{q^*}$. The variational framework for \eqref{equ} will be constructed step by step. We now prove the following propositions.
\begin{pro}\label{Pro1} Under the assumptions \eqref{n1.2},  \eqref{2.6} and ${\rm(a)}$, for any $u\in E$, there exists unique $\phi_u\in D^{1,q}(\RT)$ satisfying
\begin{eqnarray}\label{17}
  \displaystyle -\Delta_q \phi = a(x)|u|^m.
\end{eqnarray}

\end{pro}
\begin{proof} It follows from $1 < m < \frac{(q^{*} - 1)p^{*}}{q^{*}}$ that
\begin{eqnarray*}
1< \frac{q^{*}}{q^{*}-1}<\frac{mq^{*}}{q^*-1}<p^*.\label{3.2}
\end{eqnarray*}
For any $u \in E$ and $v \in D^{1,q}(\mathbb{R}^{3})$, by H\"{o}lder  and Lemma \ref{lem2.1}, we obtain
\begin{eqnarray}
    \begin{aligned}
&\left|\INT a(x)|u|^mvdx\right|\\
\leq&\left(\INT \left((a(x))^{\frac{q^*-1}{q^*}}|u|^m\right)^\frac{q*}{q*-1}dx\right)^{\frac{q^*-1}{q^*}}
\left(\INT \left((a(x))^{\frac{1}{q^*}}|v|\right)^{{q^*}} dx\right)^{\frac{1}{q^*}}\\
=&\left(\INT a(x)|u|^{\frac{mq^*}{q^*-1}}dx\right)^{\frac{q^*-1}{q^*}}\left(\INT a(x)|v|^{q^*}dx \right)^\frac{1}{q^*}\\
\le& C\|u\|^m\|v\|_{D^{1,q}(\RT)}.
\end{aligned}
\end{eqnarray}
Hence,    for each fixed $u\in E$, the linear functional
\begin{eqnarray*}
\mathcal{L}(v)=\INT a(x)|u|^mvdx,\ v\in D^{1,q}(\RT)
 \end{eqnarray*}
 is well-defined and continuous on $D^{1,q}(\RT)$.

When $q=2$,  the result follows immediately from the Lax-Milgram theorem. For $q\ne 2$, however, the Lax-Milgram theorem is not applicable. Instead, we use the Minty-Browder theorem. We show that the operator $\displaystyle -\Delta_q:D^{1,q}(\RT)\rightarrow(D^{1,q}(\RT))'$ is a continuous map  and satisfies
\begin{eqnarray}\label{16}
\langle-\Delta_qv_1-(-\Delta_qv_2),v_1-v_2\rangle>0, \ \ \ \forall \ v_1,v_2\in D^{1,q}(\RT), \ v_1\neq v_2,
\end{eqnarray}
and
\begin{eqnarray}\label{18}
\lim_{\|v\|\rightarrow\infty}\frac{\langle-\Delta_qv,v\rangle}{\|v\|}=\infty.
\end{eqnarray}

Inequality \eqref{18} holds since $q>1$. Now let $v_n\rightarrow v$ in $D^{1,q}(\RT)$. Then, $\{|\nabla v_n|^{q-2}\nabla v_n\}$ is bounded in $L^{\frac{q}{q-1}}(\RT)$, and $\nabla v_n(x)\rightarrow\nabla v(x)$ a.e. in $\RT$. By \cite[Proposition 5.4.7]{2013Willem}, we deduce that
\begin{eqnarray}\label{678}
|\nabla v_n|^{q-2}\nabla v_n\rightharpoonup|\nabla v|^{q-2}\nabla v\ \mbox{in}\ L^{\frac{q}{q-1}}(\RT).
\end{eqnarray}
Together \eqref{678} with $\||\nabla v_n|^{q-2}\nabla v_n\|_{\frac{q}{q-1}}\rightarrow\||\nabla v|^{q-2}\nabla v\|_{\frac{q}{q-1}}$, we infer that
\begin{eqnarray}\label{23}
\left(\INT||\nabla v_n|^{q-2}\nabla v_n-|\nabla v|^{q-2}\nabla v|^{\frac{q}{q-1}}dx\right)^{\frac{q-1}{q}}\rightarrow0.
\end{eqnarray}
For any $\varphi\in\ D^{1,q}(\RT)$, as $n\rightarrow\infty$,
\begin{eqnarray}\label{22}
\begin{aligned}
&\langle-\Delta_q v_n-(-\Delta_q v),\varphi\rangle\\
=&\INT (|\nabla v_n|^{q-2}\nabla v_n-|\nabla v|^{q-2}\nabla v)\nabla \varphi dx\\
\leq&\left(\INT||\nabla v_n|^{q-2}\nabla v_n-|\nabla v|^{q-2}\nabla v|^{\frac{q}{q-1}}dx\right)^{\frac{q-1}{q}}\left(\INT |\nabla \varphi|^q dx\right)^{\frac{1}{q}}\\
\leq&\left(\INT||\nabla v_n|^{q-2}\nabla v_n-|\nabla v|^{q-2}\nabla v|^{\frac{q}{q-1}}dx\right)^{\frac{q-1}{q}}\|\varphi\|_{D^{1,q}(\RT)}.
\end{aligned}
\end{eqnarray}
Combining \eqref{22} and \eqref{23}, we have
\begin{eqnarray}
\begin{aligned}
&\|-\Delta_q v_n-(-\Delta_q v)\|_{(D^{1,q}(\RT))'}\\
=&\sup\{{\langle-\Delta_q v_n+\Delta_q v,\varphi\rangle}|\ \varphi\in D^{1,q}(\RT),\|\varphi\|_{D^{1,q}(\RT)}=1\}\\
\leq&\left(\INT||\nabla v_n|^{q-2}\nabla v_n-|\nabla v|^{q-2}\nabla v|^{\frac{q}{q-1}}dx\right)^{\frac{q-1}{q}}\rightarrow0.
\end{aligned}
\end{eqnarray}
Consequently, $-\Delta_q$ is continuous on $D^{1,q}(\RT)$.

Now we introduce an elementary inequality: there exists $c_q>0$ such that for all $ x, y \in \mathbb{R}^3$,
\begin{eqnarray} \label{21} \left\{\begin{array}{ll}
\langle|x|^{q-2}x-|y|^{q-2}y, x-y\rangle_{\mathbb{R}^3} \ge c_q |x-y|^q   &  \hbox{for} \ 2\le q<3, \\
(|x|+|y|)^{2-q}\langle|x|^{q-2}x-|y|^{q-2}y, x-y\rangle_{\mathbb{R}^3} \ge c_q |x-y|^2  & \hbox{for} \ 1<q<2.
  \end{array}
 \right.
 \end{eqnarray}
 where $\left \langle \cdot , \cdot \right \rangle _{ \mathbb{R}^3} $ denotes the Euclidean inner product in $\mathbb{R}^3$. From \eqref{21}, it follows that
\begin{eqnarray}\label{20}
c_{q}\left\|v_{1}-v_{2}\right\|_{D^{1, q}(\RT)}^{q} \le\left\langle-\Delta_{q} v_{1}-\left(-\Delta_{q} v_{2}\right), v_{1}-v_{2}\right\rangle \text { for } 2 \le q<3,
 \end{eqnarray}
 and for $1<q<2$,
\begin{eqnarray} \label{19}
\begin{aligned}
   & c_q^{\frac{q}{2}} \|v_1 -v_2\|_{D^{1,q}(\RT)}^{q}\\
    =&
 c_q^{\frac{q}{2}} \INT |\nabla v_1 - \nabla v_2|^{q} dx \\
 \le &\INT (T(v_1,v_2))^{\frac{q}{2}} \left(|\nabla v_1|+|\nabla v_2|\right)^{\frac{q(2-q)}{2}} dx \\
 \le&\left(\INT T(v_1,v_2) dx\right)^{\frac{q}{2}} \left(\INT\left(|\nabla v_1|+|\nabla v_2|\right)^{q} dx \right)^{\frac{2-q}{2}}\\
  \le& \left(\left \langle -\Delta_{q}v_{1}-\left(-\Delta_{q} v_{2}\right)  ,v_1-v_2\right \rangle \right)^{\frac{q}{2}} \left(\INT\left(|\nabla v_1|+|\nabla v_2|\right)^{q} dx \right)^{\frac{2-q}{2}},
\end{aligned}
 \end{eqnarray}
 where   $$T(v_1,v_2)=\langle|\nabla v_1|^{q-2}\nabla v_1-|\nabla v_2|^{q-2}\nabla v_2, \nabla v_1-\nabla v_2\rangle_{ \mathbb{R}^3}.$$
Thus, from \eqref{20} and \eqref{19}, we derive that \eqref{16} holds. By the continuity of the map $-\Delta_q$, along with properties \eqref{16} and \eqref{18},  and since $\mathcal{L}\in (D^{1,q}(\RT))'$, the result follows from the Minty-Browder theorem. This completes the proof.
\end{proof}

Although it is difficult to provide an explicit expression for the solution $\phi_u$ of equation \eqref{17}, we can establish the following properties of $\phi_u$ using the uniqueness of the solution to \eqref{17}.

\begin{pro}\label{Pro2} For any  $u\in E$, the solution $\phi_u$ of \eqref{17}  satisfies the following properties:
\begin{enumerate}
\item[{\rm (i) }]

\begin{eqnarray*}
 \INT\left(\frac{1}{q}|\nabla \phi_u|^q-a(x)|u|^m\phi_u\right)dx=\min_{\phi\in D^{1,q}(\RT)}\left\{\INT\left(\frac{1}{q}|\nabla \phi|^q-a(x)|u|^m\phi\right)dx \right\} \text{ and } \phi_u\ge0;
\end{eqnarray*}
 \item[{\rm (ii) }] for $t>0$, $\phi_{tu}=t^{\frac{m}{q-1}} \phi_u;$
\item[{\rm (iii) }]  $\|\phi_u\|_{D^{1,q}(\RT)}\le C\|u\|^{\frac{m}{q-1}},$ where $C$ does not depend on $u ;$
     \item[{\rm (iv) }] if $u_n\rightharpoonup u$ in $E$, then $\phi_{u_n}\rightharpoonup \phi_u $ in $D^{1,q}(\RT)$ and $$\INT a(x)\phi_{u_n} |u_n|^{m-2}u_n\varphi dx\rightarrow \INT a(x)\phi_u|u|^{m-2}u\varphi dx,\ \forall \varphi\in E. $$
     \item[{\rm (v) }] if $u_n\rightarrow u$ in $E$, then $\phi_{u_n}\rightarrow \phi_u $ in $D^{1,q}(\RT).$
 \end{enumerate}
\end{pro}
\begin{proof}\ We only prove the part (iv), namely that if $u_n \rightharpoonup u$ in $E$, then $\phi_{u_n} \rightharpoonup \phi_u$ in $D^{1,q}(\mathbb{R}^3)$. The proof of other results  refer to \cite[Proposition 2.2]{2023DuSuWang}.

Define the functional
\begin{eqnarray*}
L\left(\phi_u\right)=\frac{1}{q} \INT\left|\nabla \phi_u\right|^q d x-\INT a(x)|u|^m\phi_u d x,
\end{eqnarray*}
where $u\in E$ and $\phi_u$ is the unique solution of \eqref{17}.  It is clear that  $L\in C^1(D^{1,q}(\RT), \mathbb{R})$. Let $\eta\in C_{0}^{\infty}(\mathbb{R}^3,[0,1])$  such that  $\eta|_{B_R}=1$, where $B_R=\{x\in \mathbb{R}^3:|x|\le R \}$.

The critical points of $L$ satisfy, for any $\varphi\in C^\infty_0(\RT)$,
\begin{eqnarray*}
\langle L^\prime(\phi_u),\varphi\rangle=\INT\left|\nabla \phi_u\right|^{q-2} \nabla \phi_u \nabla \varphi-\INT a(x)|u|^m \varphi d x=0.
\end{eqnarray*}
Suppose $u_n \rightharpoonup u$ in $E$. By part (iii), $\phi_{u_n}$ is bounded in $D^{1,q}(\RT)$, so there exists $w \in D^{1,q}(\RT)$, such that $\phi_{u_n} \rightharpoonup w$ in $D^{1,q}(\RT)$.
We claim that for all $\varphi \in C_0^{\infty}(\RT)$,
\begin{eqnarray}
\INT |\nabla w|^{q-2} \nabla w \nabla \varphi-\INT a(x)|u|^m \varphi d x=0.\label{00}
\end{eqnarray}
That is, $w$ is a solution of $-\Delta_q \phi=a(x)|u|^m$. By uniqueness, we conclude that  $w=\phi_u$, and hence $\phi_{u_n} \rightarrow \phi_u$.

To verify the claim, note that $L^{\prime}(w) \in(D^{1,q}(\RT))'$. For $\eta \in C_0^{\infty}\left(\mathbb{R}^3, [0,1]\right)$, we have $(\phi_{u_n}-w) \eta \rightharpoonup 0$ in $D^{1,q}(\RT)$, so
\begin{eqnarray*}
\langle L^{\prime}(w),\left(\phi_{u_n}-w\right) \eta\rangle\rightarrow 0.
\end{eqnarray*}
Since  $ L'\left(\phi_{u_n}\right)=0$ and $\phi_{u_n}$ is bounded in $D^{1,q}(\RT)$, it follows that
\begin{eqnarray*}
\langle L'\left(\phi_{u_n}\right), \left(\phi_{u_n}-w\right) \eta\rangle=0.
\end{eqnarray*}
Thus
\begin{eqnarray}\label{999}
\begin{aligned}
0=&\langle L^{\prime}(\phi_{u_n})-L^{\prime}(w), (\phi_{u_n}-w)\eta \rangle\\
=&\INT\left|\nabla \phi_{u_n}\right|^{q-2} \nabla \phi_{u_n} \nabla\left[\left(\phi_{u_n}-w\right) \eta\right]dx\displaystyle -\INT|\nabla w|^{q-2}\nabla w\nabla[(\phi_{u_n}-w) \eta] dx \\
&-\INT a(x)\left|u_n\right|^m\left(\phi_{u_n}-w\right) \eta d x+\INT a(x)|u|^m(\phi_{u_n}-w) \eta d x \\
=&\INT|\nabla \phi_{u_n}|^{q-2} \nabla \phi_{u_ n}[\nabla (\phi_{u_n}-w) \eta+(\phi_{u_n}-w) \nabla \eta]dx\\
&-\INT|\nabla w|^{q-2} \nabla w[\nabla (\phi_{u_n}-w) \eta+\left(\phi_{u_n}-w\right) \nabla \eta] d x\\
&-\INT a(x)\left|u_n\right|^m\left(\phi_{u_n}-w\right) \eta d x+\INT a(x){|u|}{ }^m(\phi_{u_n}-w\ ) \eta d x.
\end{aligned}
\end{eqnarray}
By H\"{o}lder inequality and the fact that the embedding $D^{1,q}(\RT)\hookrightarrow L^t_{\rm loc}(\RT)$ is compact for $1\leq t<q^*$,  we obtain
\begin{eqnarray}\label{888}
\left.
  \begin{array}{ll}
&\displaystyle \INT\left|\nabla \phi_{u_n}\right|^{q-2} \nabla \phi_{u_ n}\left(\phi_{u_n}-w\right) \nabla \eta dx\rightarrow0,\\ [1em]
&\displaystyle  \INT|\nabla w|^{q-2} \nabla w(\phi_{u_n}-w) \nabla \eta dx\rightarrow0.
  \end{array}
\right.
\end{eqnarray}
We claim that there exists  $1<t_0<q^*$  satisfying $mt'_0=m\frac{t_0}{t_0-1}\leq p^*$. Indeed,
\begin{eqnarray*}
m\frac{t_0}{t_0-1}\leq p^*\Leftrightarrow m\leq \frac{p^*(t_0-1)}{t_0}.
\end{eqnarray*}
By assumption that $m<\frac{p^*(q^*-1)}{q^*}$, the existence of such a $t_0$ is guaranteed.
Thus
\begin{eqnarray}\label{777}
\begin{aligned}
&\left|\INT a(x)\left|u_n\right|^m\left(\phi_{u_n}-w\right) \eta dx\right|\\
\le&\left(\INT a(x)|u|^{m t_0^{\prime}}\eta  dx\right)^{\frac{1}{t'_0}}\left(\INT a(x)|\phi_{u_n}-w|^{t_0} \eta dx\right)^ \frac{1}{t_0}\rightarrow 0.
\end{aligned}
\end{eqnarray}
Since $\phi_{u_n} \rightharpoonup w $ in $L^{q^*} $ and $\INT \left(a(x)|u|^m \eta\right)^{\frac{q^*}{q^*-1}} dx<\infty$, we get that
\begin{eqnarray}\label{666}
\INT a(x){|u|}{ }^m(\phi_{u_n}-w\ ) \eta d x\rightarrow0.
\end{eqnarray}
By \eqref{888}-\eqref{666} and \eqref{999}, we get
\begin{eqnarray*}
\begin{aligned}
0=&\INT\left|\nabla \phi_{u_n}\right|^{q-2} \nabla \phi_{u_n} {\nabla}\left(\phi_{u_n}-w\right) \eta dx-\INT|\nabla w|^{q-2} \nabla w\nabla\left(\phi_{u_n}-w\right) \eta dx\\
=&\INT\left(\left| \nabla \phi_{u_n}\right|^{q-2} \nabla \phi_{u_n}-|\nabla w|^{q-2} \nabla w\right)[\nabla (\phi_{u_n}-w)]\eta d x.
\end{aligned}
\end{eqnarray*}
By \eqref{21} and the definition of $\eta$, we obtain that
\begin{eqnarray*}
    \int_{B_R} (|\nabla\phi_{u_n}|^{q-2}\nabla\phi_{u_n}-|\nabla w|^{q-2}\nabla w)(\nabla \phi_{u_n}-\nabla w)dx  \rightarrow 0.
\end{eqnarray*}
Using \eqref{21} and arguing as in \eqref{20} and \eqref{19}, we conclude
\begin{eqnarray*}
    \int_{B_R}|\nabla \phi_{u_n}-\nabla w|^qdx\rightarrow0.
\end{eqnarray*}
By the  arbitrariness of $B_R$, we deduce that  $ \nabla \phi_{u_n}(x) \rightarrow \nabla w(x)$ a.e. in  $\mathbb{R}^3$.
For any   $\varphi \in C_0^{\infty}(\RT)$, we have
\begin{eqnarray}
0=\langle L^{\prime}\left(\phi_{u_n}\right), \varphi\rangle=\INT\left|\nabla \phi_{u_n}\right|^{q-2} \nabla \phi_{u_n}\nabla \varphi dx-\INT a(x)\left|u_n\right|^m \varphi d x.\label{nn3.16}
\end{eqnarray}
By $u_n\rightharpoonup u$ in $E$ and \eqref{n2.1}, we have
\begin{eqnarray*}
u_n\rightharpoonup u\ {\rm in}\ D^{1,p}(\RT).
\end{eqnarray*}
As $1<m<\frac{(q^*-1)p^*}{q^*}$, it follows that
\begin{eqnarray}
 u_n\to u\ {\rm in}\ L^m_{\rm loc}(\RT).\label{n3.16}
\end{eqnarray}
Combining with (a), \eqref{n3.16}, $\varphi \in C_0^{\infty}(\RT)$ and the Lebesgue's dominated convergence theorem, we obtain,
\begin{eqnarray*}
\INT a(x)\left|u_n\right|^m \varphi d x\to \INT a(x)\left|u\right|^m \varphi d x.
\end{eqnarray*}
Since $\nabla \phi_{u_n} (x)\rightarrow \nabla w(x)$ a.e. in  $\mathbb{R}^3$,  it follows from \eqref{nn3.16} that
$$
0=\INT|\nabla w|^{q-2} \nabla w \nabla \eta dx-\INT a(x)|u|^m \eta d x.
$$
This establishes \eqref{00}, and hence $\phi_{u_n}\rightharpoonup \phi_u $ in $D^{1,q}(\RT)$. The proof is complete.

\end{proof}

We define the functional
\begin{eqnarray*}
\begin{aligned}
  I(u)=&\frac{1}{p}\INT \left | \nabla u  \right | ^pdx-\frac{\lambda }{p}\INT g(x)\left | u \right |^pdx  +\frac{q-1}{qm}\INT a(x)\phi_{u} \left |  u\right | ^{m} dx\\
&-\frac{1 }{r} \INT k(x)  \left |u  \right |^{r}dx +\frac{1}{s}\INT h(x)\left | u \right |  ^{s}dx.
\end{aligned}
\end{eqnarray*}
We will prove that $I\in C^1(E,\mathbb{R})$ and  that the critical points of $I$
correspond to solutions of \eqref{equ}. The key step is to verify  the map $u \longmapsto \phi_u \in C^1(E,D^{1,q}(\RT))$. The ideas for the following results  are adapted from \cite[Proposition 2.3]{2023DuSuWang}.

\begin{pro}\label{pro3.4}
$I\in C^1(E,\mathbb{R})$, and for any $u,v\in E$,
\begin{eqnarray*}
\begin{aligned}
    \langle I'(u),v\rangle=&\INT\left|\nabla u \right|^{p-2}\nabla u\nabla vdx-\lambda\INT g(x)\left | u \right |^{p-2}uvdx  +\INT a(x)\phi_{u} \left |  u\right | ^{m-2}uvdx\\
  &- \INT k(x)\left|u\right|^{r-2} uvdx +\INT h(x)\left | u \right |  ^{s-2}uvdx.
  \end{aligned}
\end{eqnarray*}
\end{pro}
\begin{proof} We first prove that $I$ is G$\hat{a}$teaux differentiable, that is
 \begin{eqnarray}\label{12}
\lim_{t\rightarrow0}\frac{ {I}(u+tv)-I(u)}{t}=\langle I'_{G}(u),v\rangle,\ \quad u,v\in E,
 \end{eqnarray}
where
\begin{eqnarray*}
\begin{aligned}
  \langle I'_{G}\left(u\right),v\rangle=&\INT\left|\nabla u \right|^{p-2}\nabla u\nabla  vdx-\lambda\INT g(x)\left| u \right|^{p-2} uvdx  +\INT a(x)\phi_{u} \left| u\right|^{m-2} uvdx\\
  &-\INT k(x)\left|u\right|^{r-2}uvdx +\INT h(x)\left|u \right|^{s-2}uvdx.
  \end{aligned}
\end{eqnarray*}
Clearly,  $\langle I'_{G}(u), v\rangle$ is a continuous linear functional  in $v$. Define
$$I(u+tv)-I(u)-t\langle I'_{G}(u), v\rangle =\mathcal{A}+\mathcal{B}+\mathcal{C}+\mathcal{D},$$
where
 \begin{eqnarray*}
 \begin{aligned}
    \mathcal{A}:\triangleq&\frac{1}{p}\INT (|\nabla (u+tv)|^p-|\nabla u|^p)dx-\frac{\lambda}{p}\INT g(x)(|u+ tv|^p-|u|^p)dx-t\INT |\nabla u|^{p-2}\nabla u\nabla vdx\\
&+t\lambda\INT g(x)|u|^{p-2} uvdx,
 \end{aligned}
\end{eqnarray*}
\begin{flalign*}
&\mathcal{B}:\triangleq-\frac{1}{r}\INT k(x)|u+tv|^rdx+\frac{1}{r}\INT k(x)|u|^r dx+t\INT k(x) |u|^{r-2}uv dx,&
\end{flalign*}

\begin{flalign*}
& \mathcal{C}:\triangleq-\frac{1}{s}\INT h(x)|u+tv|^sdx+\frac{1}{s}\INT h(x)|u|^s dx+t\INT h(x) |u|^{s-2} uvdx,&
\end{flalign*}

\begin{flalign*}
&\mathcal{D}:\triangleq\frac{q-1}{qm} \INT a(x)\phi_{u+tv}|u+tv|^mdx-\frac{q-1}{qm} \INT a(x)\phi_u|u|^mdx-t\INT a(x)\phi_u|u|^{m-2}uvdx.&
\end{flalign*}

It is straightforward to verify that
$$\mathcal{A}=o(t),\mathcal{B}=o(t),\mathcal{C}=o(t),\quad t\rightarrow 0.$$
To prove \eqref{12},  it suffices to show that
\begin{eqnarray}\label{28}
    \mathcal{D}=o(t), \quad t\rightarrow 0.
\end{eqnarray}
From Proposition  \ref{Pro2} (i), we have
\begin{eqnarray}\label{3.9}
\mathcal{M}(u,\phi_u)=\min_{\phi\in\D}\mathcal{M} (u,\phi),
\end{eqnarray}
where
$$\mathcal{M}(u,\phi)=\frac{1}{q}\INT |\nabla \phi|^q dx-\INT a(x)\phi |u|^m dx.$$
Using $\mathcal{M}(u+tv,\phi_u) \ge \mathcal{M}(u+tv,\phi_{u+tv})$, we deduce that
\begin{eqnarray}
\begin{aligned}
    \mathcal{D}&\displaystyle =\frac{q-1}{qm} \INT a(x)\phi_{u+tv}|u+tv|^mdx-\frac{q-1}{qm} \INT a(x)\phi_u|u|^mdx-t\INT a(x)\phi_u|u|^{m-2} uvdx \\
&\displaystyle =-\frac{1}{m}\mathcal{M}(u+tv,\phi_{u+tv})-\frac{q-1}{qm} \INT a(x)\phi_u|u|^mdx-t\INT a(x)\phi_u|u|^{m-2}uvdx \\
&\displaystyle \geq-\frac{1}{m}\mathcal{M}(u+tv,\phi_u)-\frac{q-1}{qm} \INT a(x)\phi_u|u|^mdx-t\INT a(x)\phi_u|u|^{m-2}uvdx \\
&\displaystyle  =\frac{1}{m}\INT a(x)\phi_u(|u+tv|^m-|u|^m-mt|u|^{m-2}uv)dx.
\label{8}
\end{aligned}
\end{eqnarray}
We claim that
\begin{eqnarray}\label{7}
  \INT a(x)\phi_u(|u+tv|^m-|u|^m-mt|u|^{m-2}uv)dx=o(t),\quad t\rightarrow 0.
\end{eqnarray}
Indeed, since $m>1$,  for a.e. $x\in \mathbb{R}^3$,
\begin{eqnarray}\label{9}
\lim_{t\rightarrow0}\frac{a(x)\phi_u(x)(|u(x)+tv(x)|^m-|u(x)|^m-mt|u(x)|^{m-2}u(x)v(x))}{t}=0.
\end{eqnarray}
Using the elementary inequality: for any $\tau>0$, there exists $C_\tau>0$,
such that for any $a,b\in \mathbb{R}$
\begin{eqnarray}  \label{11}
|a+b|^\tau\leq C_\tau(|a|^\tau+|b|^\tau).
\end{eqnarray}
Then, as $m>1$, by the Lagrange theorem, there exists $\theta\in \mathbb{R}$ such that $|\theta|\le|t|$ and
\begin{eqnarray}\label{5}
\begin{aligned}
& \frac{1}{t}\left|a(x)\phi_{u}(x)\left[|u(x)+t v(x)|^{m}-|u(x)|^{m}-t m|u(x)|^{m-2} u(x) v(x)\right]\right| \\
\le & m a(x)\phi_{u}(x)(||u(x)+\left.\theta v(x)|^{m-2}(u(x)+\theta v(x)) v(x)|+|| u(x)|^{m-2} u(x) v(x)|\right)  \\
\leqslant & C a(x)\phi_{u}(x)\left(|u(x)|^{m-1}|v(x)|+|v(x)|^{m}\right).
\end{aligned}
\end{eqnarray}
By  H\"{o}lder and Sobolev inequalities, using \eqref{11} and Proposition \ref{Pro2} (iii), we have
\begin{eqnarray}\label{10}
\begin{aligned}
&\left|\INT a(x)\phi_u(|u|^{m-1}|v|+|v|^m)dx\right|\\
\leq&\left(\INT |\phi_u|^{q^*}dx\right)^{\frac{1}{q^*}}
\left(\INT (a(x))^{\frac{q^*}{q^*-1}} (|u|^{m-1}|v|+|v|^m)^{\frac{q^*}{q^*-1}} dx\right)^{\frac{q^*-1}{q^*}}\\
\leq& C\|\phi_u\|_{D^{1,q}(\RT)}\left(\INT (a(x))^{\frac{q^*}{q^*-1}} (|u|^{m-1}|v|+|v|^m)^{\frac{q^*}{q^*-1}} dx\right)^{\frac{q^*-1}{q^*}}\\
\leq& C\|u\|^{\frac{m}{q-1}}
\left(\INT (a(x)|u|^{m-1}|v|)^{\frac{q^*}{q^*-1}}+(a(x))^{\frac{q^*}{q^*-1}}|v|^{\frac{mq^*}{q^*-1}} dx\right)^{\frac{q^*-1}{q^*}}.
\end{aligned}
\end{eqnarray}
Since $1<m<\frac{(q^*-1)p^*}{q^*}$, by  H\"{o}lder inequality and  Lemma \ref{lem2.1}, we have
\begin{eqnarray}\label{1}
\begin{aligned}
&\INT |a(x)|u|^{m-1}v|^{\frac{q^*}{q^*-1}} dx\\
\leq&\|a\|_\infty^{\frac{1}{q^*-1}}\INT a(x)\left||u|^{m-1}v\right|^{\frac{q^*}{q^*-1}} dx\\
\leq&\|a\|_\infty^{\frac{1}{q^*-1}}\left(\INT a(x)|u|^{\frac{mq^*}{q^*-1}} dx\right)^{\frac{m-1}{m}}\left(\INT a(x)|v|^{\frac{mq^*}{q^*-1}} dx\right)^{\frac{1}{m}}\\
<&\infty.
\end{aligned}
\end{eqnarray}
Similarly,
\begin{eqnarray}
\INT (a(x))^{\frac{q^*}{q^*-1}}|v|^{\frac{mq^*}{q^*-1}} dx<\infty.\label{n3.28}
\end{eqnarray}
From \eqref{10}, \eqref{1}  and \eqref{n3.28}, the function $a(x)\phi_u(|u|^{m-1}|v|+|v|^m)$ is in $L^1(\RT)$. Together this with \eqref{5} and \eqref{9},  the dominated convergence theorem implies \eqref{7}.  From \eqref{8} and \eqref{7}, we deduce that
\begin{eqnarray}\label{27}
    \mathcal{D}\ge o(t).
\end{eqnarray}
We now prove that
    \begin{eqnarray}\label{26}
    \mathcal{D}\le o(t).
\end{eqnarray}
Actually, using $\mathcal{M}(u,\phi_u)\leq \mathcal{M}(u,\phi_{u+tv})$, by  H\"{o}lder inequality, we deduce that
\begin{eqnarray}\label{25}
\begin{aligned}
\mathcal{D}=&\frac{q-1}{mq} \INT a(x)\phi_{u+tv}|u+tv|^mdx-\frac{q-1}{mq} \INT a(x)\phi_u|u|^mdx-t\INT a(x)\phi_u|u|^{m-2}uvdx\\
=&\frac{q-1}{mq} \INT a(x)\phi_{u+tv}|u+tv|^mdx+\frac{1}{m}\mathcal{M}(u,\phi_u)-t\INT a(x)\phi_u|u|^{m-2}uvdx\\
\leq&\frac{q-1}{mq} \INT a(x)\phi_{u+tv}|u+tv|^mdx+\frac{1}{m}\mathcal{M}(u,\phi_{u+tv})-t\INT a(x)\phi_u|u|^{m-2}uvdx\\
=&\frac{q-1}{mq} \INT a(x)\phi_{u+tv}|u+tv|^mdx+\frac{1}{m}\left\{\frac{1}{q}\INT |\nabla \phi_{u+tv}|^q dx-\INT a(x)|u|^m\phi_{u+tv}dx\right\}\\
&-t\INT a(x)\phi_u|u|^{m-2}uvdx\\
=&\frac{1}{m}\INT [a(x)\phi_{u+tv}(|u+tv|^m-|u|^m)-mta(x)\phi_u|u|^{m-2}uv]dx\\
=&\frac{1}{m}\INT [a(x)\phi_{u+tv}(|u+tv|^m-|u|^m-mt|u|^{m-2}uv)+mta(x)(\phi_{u+tv}-\phi_u)|u|^{m-2}uv]dx\\
\leq&\frac{1}{m}\left(\INT |\phi_{u+tv}|^{q^*}dx\right)^{\frac{1}{q^*}}
\left(\INT (a(x))^{\frac{q^*}{q^*-1}}||u+tv|^m-|u|^m-mt|u|^{m-2}uv|^{\frac{q^*}{q^*-1}} dx\right)^{\frac{q^*-1}{q^*}}\\
&+ t\left(\INT |\phi_{u+tv}-\phi_u|^{q^*}dx\right)^{\frac{1}{q^*}}
\left(\INT (a(x))^{\frac{q^*}{q^*-1}}||u|^{m-2}uv|^{\frac{q^*}{q^*-1}} dx\right)^{\frac{q^*-1}{q^*}}.
\end{aligned}
\end{eqnarray}
We first verify that
\begin{eqnarray} \label{6}
\lim_{t\rightarrow0}\frac{1}{t}\left(\INT (a(x))^{\frac{q^*}{q^*-1}}||u+tv|^m-|u|^m-mt|u|^{m-2}uv|^{\frac{q^*}{q^*-1}} dx\right)^{\frac{q^*-1}{q^*}}=0.
\end{eqnarray}
To prove \eqref{6}, it suffices to show that
\begin{eqnarray} \label{4}
\lim_{t\rightarrow0}t^{-\frac{q^*}{q^*-1}}\INT (a(x))^{\frac{q^*}{q^*-1}}||u+tv|^m-|u|^m-mt|u|^{m-2}uv|^{\frac{q^*}{q^*-1}} dx=0.
\end{eqnarray}
Indeed, since $m>1$, it is easy to see that for a.e. $x\in\mathbb{R}^3$,
\begin{eqnarray}\label{3}
   \lim_{t\rightarrow0}\frac{a(x)(|u(x)+tv(x)|^m-|u(x)|^m-mt|u(x)|^{m-2}u(x)v(x))}{t}=0.
\end{eqnarray}
As in \eqref{5}, by the Lagrange theorem, there exists $\theta\in \mathbb{R} $ such that $|\theta|\leq|t|$ and
\begin{eqnarray}\label{2}
\begin{aligned}
&\frac{1}{t}a(x)\left||u(x)+t v(x)|^{m}-|u(x)|^{m}-mt| u(x)|^{m-2} u(x) v(x) \right\rvert\\
\le& Ca(x)\left(|u(x)|^{m-1}|v(x)|+|v(x)|^{m}\right).
\end{aligned}
\end{eqnarray}
From \eqref{1}, we infer that the function $(a(x)|u|^{m-1}|v|+a(x)|v|^m)^{\frac{q^*}{q^*-1}}\ \mbox{is \ in}\ L^1(\RT).$ Combining this with \eqref{2} and \eqref{3},  the dominated convergence theorem implies \eqref{4}. Applying Proposition \ref{Pro2} (v) and the Sobolev embedding theorem, we   derive that
\begin{eqnarray}\label{24}
    \lim_{t\rightarrow0}\INT|\phi_{u+tv}-\phi_u|^{q^*}dx=0.
\end{eqnarray}
Using \eqref{6}  and  \eqref{24}, it follows from  \eqref{25} that  \eqref{26}. Then \eqref{26} and \eqref{27} yields \eqref{28}.

To end the proof, we show the continuity of $I'_{G}:E\rightarrow (E)'$.  It is sufficient to prove that
$$\langle\Psi'_{G}(u),v\rangle:=\INT a(x)\phi_{u}|u|^{m-2}uvdx$$
is continuous on $E$. Let $u_n\rightarrow u$ in $E$. By Proposition \ref{Pro2} (v) and the Sobolev embedding theorem, we infer that
\begin{eqnarray}
\phi_{u_n} \rightarrow \phi_u\ {\rm in}\ L^{q^*}(\RT). \label{n3.37}
\end{eqnarray}
Moreover, based on Lemma \ref{lem2.1} and  the dominated convergence theorem,
\begin{eqnarray}
a(x)|u_n|^{m-2}u_n-a(x)|u|^{m-2}u\rightarrow0 \ \mbox{in}\ L^{\frac{mq^*}{(m-1)(q^*-1)}}(\RT).\label{n3.38}
\end{eqnarray}
Then, for any $v\in E$, by H\"{o}lder inequality,   \eqref{n3.37} and \eqref{n3.38}, we obtain that
\begin{eqnarray*}
\begin{aligned}
&\langle\Psi'_{G}(u_n)-\Psi'_{G}(u),v\rangle\\
    =&\INT a(x)(\phi_{u_n}|u_n|^{m-2}u_n-\phi_{u}|u|^{m-2}u)vdx\\
    =&\INT\phi_{u_n}(a(x)|u_n|^{m-2}u_n-a(x)|u|^{m-2}u    )vdx+\INT a(x)(\phi_{u_n}-\phi_{u} )|u|^{m-2}uvdx\\
    \le&\left(\INT |\phi_{u_n}|^{q^*}dx\right)^{\frac{1}{q^*}}\left(\INT a(x)\left||u_n|^{m-2}u_n-|u|^{m-2}u\right|^{\frac{mq^*}{(m-1)(q^*-1)}}dx\right)^{\frac{(m-1)(q^*-1)}{mq^*}}\\
    &\left(\INT (a(x))^{\frac{q^*+m-1}{q^*-1}}|v|^{\frac{mq^*}{q^*-1}}dx\right)^{\frac{q^*-1}{mq^*}}\\
    &+\left(\INT |\phi_{u_n}-\phi_u|^{q^*}dx\right)^{\frac{1}{q^*}}\left(\INT a(x)||u|^{m-2}u|^{\frac{mq^*}{(m-1)(q^*-1)}}dx\right)^{\frac{(m-1)(q^*-1)}{mq^*}}\\
    &\left(\INT (a(x))^{\frac{q^*+m-1}{q^*-1}}|v|^{\frac{mq^*}{q^*-1}}dx\right)^{\frac{q^*-1}{mq^*}}\rightarrow 0.
    \end{aligned}
    \end{eqnarray*}
Thus,
\begin{eqnarray*}
\begin{aligned}
&\|\Psi_G'(u_n)-\Psi_G'(u)\|_{(E)'}\\
=&\sup\{\langle\Phi_G'(u_n)-\Phi_G'(u),v\rangle|\ v\in E,\|v\|=1\}\\
\leq&C\left(\INT a(x)||u_n|^{m-2}u_n-|u|^{m-2}u|^{\frac{mq^*}{(m-1)(q^*-1)}}dx\right)^{\frac{(m-1)(q^*-1)}{mq^*}}\left(\INT |\phi_{u_n}|^{q^*}dx\right)^{\frac{1}{q^*}}\\
&+C\left(\INT |\phi_{u_n}-\phi_u|^{q^*}dx\right)^{\frac{1}{q^*}}\left(\INT a(x)||u|^{m-2}u|^{\frac{mq^*}{(m-1)(q^*-1)}}dx\right)^{\frac{(m-1)(q^*-1)}{mq^*}}\\
\rightarrow&0.
\end{aligned}
\end{eqnarray*}
This completes the proof.
\end{proof}

\begin{pro}\label{pro3.2}
    Let $\{u_n\}\subset E$ and $u_n\rightharpoonup u$ in $E$ such that $I'(u_n)\rightarrow0$. Then, there exists $u\in E$ so that, up to a subsequence $\nabla u_n(x)\rightarrow \nabla u(x)$ a.e. $x\in \mathbb{R}^3$.
\end{pro}
\begin{proof}
Since $u_n\rightharpoonup u$ in $E$, we may extract a subsequence (still denoted by $\{u_n\}$) such that
\begin{eqnarray}\label{3311}
   \left\{\begin{array}{ll}
u_n\rightharpoonup u & \text { in } D^{1,p}(\RT), \\
u_{n} \rightarrow u & \text { in } L_{l}^{\gamma }(\RT), \forall 1 < \gamma <p^{*}, \\
u_{n} \rightarrow u & \text { in } L_{\rm loc}^q(\RT),\forall 1\le q<p^*, \\
u_{n}(x) \rightarrow u(x) & \text { a.e. } x \in \mathbb{R}^{3} .
\end{array}\right.
\end{eqnarray}
We adapt some ideas from \cite{1992BoccardoMurat}. Let $\eta\in C_{0}^{\infty}(\mathbb{R}^3,[0,1])$ such that $\eta|_{B_R}=1$ and supp $\eta\subset B_{2R}$, where $B_R=\{x\in \mathbb{R}^3:|x|\le R \}$. Since $I'(u)\in (E)'$, $(u_n-u)\eta\rightharpoonup 0 $ in $E$, and $I'(u_n)\rightarrow 0$, we have
\begin{eqnarray}\label{3312}
\begin{aligned}
o(1)=&\langle I'(u_n)-I'(u),(u_n-u)\eta \rangle\\
=&\INT(|\nabla u_n|^{p-2}\nabla u_n-|\nabla u|^{p-2}\nabla u)\nabla
[(u_n-u)\eta]dx \\
&-\lambda\INT g(x)(|u_n|^{p-2} u_n-| u|^{p-2}u)
(u_n-u)\eta dx\\
&-\INT k(x)(|u_n|^{r-2} u_n-| u|^{r-2}u)
(u_n-u)\eta dx\\
&+\INT a(x)(\phi_{u_n}|u_n|^{m-2} u_n-\phi_{u}| u|^{m-2}u)
(u_n-u)\eta dx\\
&-\INT h(x)(|u_n|^{s-2} u_n-| u|^{s-2}u)
(u_n-u)\eta dx.
\end{aligned}
\end{eqnarray}
By H\"{o}lder inequality, using \eqref{3311},  we obtain
\begin{eqnarray}\label{3313}
\begin{aligned}
 \lambda\INT g(x)(|u_n|^{p-2} u_n-| u|^{p-2}u)
(u_n-u)\eta dx=o(1),\\
\INT k(x)(|u_n|^{r-2} u_n-| u|^{r-2}u)
(u_n-u)\eta dx=o(1),\\
\INT h(x)(|u_n|^{s-2} u_n-| u|^{s-2}u)
(u_n-u)\eta dx=o(1),
\end{aligned}
\end{eqnarray}
and
\begin{eqnarray}\label{3314}
\begin{aligned}
&\INT|\nabla u_n|^{p-2}\nabla u_n\nabla \eta
(u_n-u)dx\\
\le& \left(\int_{B_{2R}}\left||\nabla u_n^{}|^{p-2}\nabla u\nabla\eta \right|^{\frac{p}{p-1}}dx \right)^{\frac{p-1}{p}}\left(\int_{B_{2R}}|u_n-u|^pdx  \right)^{\frac{1}{p}}\\
\le&\left(\int_{B_{2R}}|\nabla u_n|^pdx \right)^{\frac{p-1}{p}}\left(\int_{B_{2R}}|u_n-u|^pdx  \right)^{\frac{1}{p}}\\
=&o(1).
\end{aligned}
\end{eqnarray}
Moreover, by H\"{o}lder inequality, we have
\begin{eqnarray*}
\begin{aligned}
&\INT a(x)(\phi_{u_n}|u_n|^{m-2} u_n-\phi_{u}| u|^{m-2}u)
(u_n-u)\eta dx\nonumber\\
\le&\left(\INT |\phi_{u_n}|^{q^*} \right)^{\frac{1}{q^*}}\left[\int_{B_{2R}}\left( a(x)\left||u_n|^{m-2}u_n(u_n-u)\eta\right| \right)^{\frac{q^*}{q^*-1}} \right]^{\frac{q^*-1}{q^*}}\nonumber\\
&+\left(\INT |\phi_{u}|^{q^*} \right)^{\frac{1}{q^*}}\left[\int_{B_{2R}}\left( a(x)\left||u|^{m-2}u(u_n-u)\eta\right| \right)^{\frac{q^*}{q^*-1}} \right]^{\frac{q^*-1}{q^*}}\nonumber\\
\le&\|\phi_{u_n}\|_{q^*}\left[\int_{B_{2R}} a(x)(u_n-u)^{\frac{mq^*}{q^*-1}} dx\right]^{\frac{q^*-1}{mq^*}}\left[\int_{B_{2R}} (a(x))^{1+\frac{m}{(q^*-1)(m-1)}} |u_n|^{\frac{mq^*}{q^*-1}}dx \right]^{\frac{(q^*-1)(m-1)}{mq^*}}\nonumber\\
&+\|\phi_{u}\|_{q^*}\left[\int_{B_{2R}} a(x)(u_n-u)^{\frac{mq^*}{q^*-1}} dx\right]^{\frac{q^*-1}{mq^*}}\left[\int_{B_{2R}} (a(x))^{1+\frac{m}{(q^*-1)(m-1)}} |u|^{\frac{mq^*}{q^*-1}}dx \right]^{\frac{(q^*-1)(m-1)}{mq^*}}.
\end{aligned}
\end{eqnarray*}
Thus, we deduce that
\begin{eqnarray}\label{3315}
    \INT a(x)(\phi_{u_n}|u_n|^{m-2} u_n-\phi_{u}| u|^{m-2}u)
(u_n-u)\eta dx=o(1).
\end{eqnarray}
From \eqref{3312}-\eqref{3315}, it follows that
\begin{eqnarray*}
\INT(|\nabla u_n|^{p-2}\nabla u_n-|\nabla u|^{p-2}\nabla u)\nabla
[(u_n-u)\eta]dx=o(1).
\end{eqnarray*}
Using again \eqref{21} and arguing as in \eqref{20} and \eqref{19}, we deduce that
\begin{eqnarray*}
    \lim_{n\rightarrow\infty}\int_{B_R}|\nabla u_n-\nabla u|^pdx=0.
\end{eqnarray*}
Thus, up to a subsequence,
$$\nabla u_n(x)\rightarrow \nabla u(x) \text{ a.e. } x\in B_R.$$
In view of the arbitrariness of $B_R$, we conclude that
$$\nabla u_n(x)\rightarrow \nabla u(x) \text{ a.e. } x\in \mathbb{R}^3.$$
The proof is complete.
\end{proof}

\begin{pro}
    Let $(u,\phi)\in E\times D^{1,q}(\RT)$. Then, $(u,\phi)$ is a solution of \eqref{equ} if and only if $u$ is a critical point of $I$ and $\phi=\phi_u$.
\end{pro}
\begin{proof} Define the functional $\mathcal{H}: E \times D^{1,q}(\mathbb{R}^3) \to \mathbb{R}$ by
\begin{eqnarray*}
\begin{aligned}
\mathcal{H}(u,\phi)=&\frac{1}{p}\INT \left | \nabla u  \right | ^p dx-\frac{\lambda }{p}\INT g(x)\left | u \right |^p dx  +\frac{1}{m}\INT a(x)\phi_{u} \left |  u\right | ^{m} d x-\frac{1}{q m}\INT|\nabla \phi|^q dx\\
&-\frac{1 }{r} \INT k(x)  \left |u  \right |^{r}d x +\frac{1}{s}\INT h(x)\left | u \right |  ^{s}dx.
\end{aligned}
\end{eqnarray*}
For any $v,w\in E\times D^{1,q}(\RT)$,
\begin{eqnarray*}
\begin{aligned}
    \langle\partial_u \mathcal{H}(u,\phi),v\rangle=&\INT\left|\nabla u \right|^{p-2}\nabla u\nabla  vdx-\lambda\INT g(x)\left | u \right |^{p-2} uvdx  +\INT a(x)\phi|u|^{m-2}uvdx\\
&- \INT k(x)\left|u\right|^{r-2}uv dx +\INT h(x)\left|u \right|^{s-2}uv dx,
\end{aligned}
\end{eqnarray*}
\begin{eqnarray*}
  \langle\partial_\phi \mathcal{H}(u,\phi),w\rangle=\frac{1}{m}\INT a(x)|u|^mw dx-\frac{1}{m}\INT |\nabla\phi|^{q-2}\nabla\phi\nabla w dx.
\end{eqnarray*}
Therefore,
\begin{center}
$(u,\phi)$ is a solution of \eqref{equ} $\Leftrightarrow \partial_u \mathcal{H}(u,\phi)=0$ and $\partial_\phi \mathcal{H}(u,\phi)=0$.\\
\end{center}
From the definition of $I'(u)$ given in Proposition \ref{pro3.4}, it follows that
\begin{center}
$\partial_u \mathcal{H}(u,\phi)=0 $ and $ \partial_\phi \mathcal{H}(u,\phi)=0 \Leftrightarrow I'(u)=0$ and $\phi=\phi_u$.
\end{center} \end{proof}

\section{Proof of Theorem 1.1}\label{sec4}
We first show that the functional $I$ is coercive and satisfies the Palais-Smale condition.
\begin{lem}\label{lem4.1}
    Under the assumption of Theorem \ref{thm1.1}, the functional $I$ is coercive.
\end{lem}
\begin{proof}
By H\"{o}lder and Sobolev inequalities, using \eqref{38}, we obtain
\begin{eqnarray}\label{41}
\begin{aligned}
   \INT k(x)|u|^rdx
    &\le \left(\INT |k(x)|^{r_0}dx\right)^{\frac{1}{r_0}}\left( \INT |u|^{r\frac{r_0}{r_0-1}}dx  \right)^{\frac{r_0-1}{r_0}}\\
    &=\|k\|_{r_0}\|u\|^r_{p^*}\\
    &\le\|k\|_{r_0}S_{p}^{-{\frac{r}{p}}}\|u\|^r_{D^{1,p}(\RT)},
    \end{aligned}
\end{eqnarray}
where $r_0$ is defined in \eqref{456}. If $0\le\lambda<\lambda_1$, then by Lemma \ref{lem2.2} (ii), there exists $c > 0$ such that for $\|u\|>1$,
\begin{eqnarray}
  \begin{aligned}
    I(u)&\displaystyle =\frac{1}{p}\INT \left| \nabla u  \right| ^pdx-\frac{\lambda }{p}\INT g(x)\left| u \right|^pdx  +\frac{q-1}{qm}\INT a(x)\phi_{u} \left|  u\right| ^{m} dx\\
&\displaystyle \ \ \ -\frac{1 }{r} \INT k(x)  \left|u  \right|^{r}dx +\frac{1}{s}\INT h(x)\left| u \right|  ^{s}dx\\
&\displaystyle \ge \frac{1}{p}\left(1-\frac{\lambda}{\lambda_1}\right)\INT|\nabla u|^pdx+\frac{1}{s}\INT h(x)\left|u\right|^{s}dx\\
&\displaystyle\ \ \  -\frac{1 }{r} \INT k(x)  \left|u  \right|^{r}dx+\frac{q-1}{qm}\INT a(x)\phi_{u}\left|u\right|^{m} dx\\
&\displaystyle \ge c\|u\|^p-\frac{1}{r}\|k\|_{r_0}S_{p}^{-{\frac{r}{p}}}\|u\|^r_{D^{1,p}}+\frac{q-1}{qm}\INT a(x)\phi_{u} \left|  u\right|^{m} dx\\
&\displaystyle \ge c\|u\|^p-\frac{1}{r}\|k\|_{r_0}S_{p}^{-{\frac{r}{p}}}\|u\|^r+\frac{q-1}{qm}\INT a(x)\phi_{u} \left|  u\right|^{m} dx\\
&\displaystyle \ge c\|u\|^p-\frac{1}{r}\|k\|_{r_0}S_{p}^{-{\frac{r}{p}}}\|u\|^r.
    \end{aligned}
\end{eqnarray}
Since $p>r$,  it follows that $I$ is coercive. A symmetric argument holds for $-\lambda_1<\lambda<0$. \end{proof}

 In general, to prove the (PS) condition, the reflexivity of the working space is typically required. Although it is unclear whether $E$ is reflexive, we can still establish the following result.

\begin{lem}\label{lem4.2}
  Under the assumptions of Theorem \ref{thm1.1}, the functional $I$ satisfies the {\rm(PS)} condition.
\end{lem}
\begin{proof}
By Lemma \ref{lem4.1}, every (PS) sequence $ \{u_n\} $ of $I$ is bounded in $E$. Thus $\{u_n\}$ is bounded in $D^{1,p}(\RT)$  and $\{(h(x))^{\frac{1}{s}}u_n\}$ is bounded in $L^s(\RT)$. Hence, up to a subsequence, we have
\begin{center}
  $u_n\rightharpoonup u$ in $D^{1,p}(\RT)$ , $(h(x))^{\frac{1}{s}}u_n\rightharpoonup (h(x))^{\frac{1}{s}}u $ in $L^s(\RT)$.
\end{center}
Consequently,
\begin{eqnarray}\label{43}
    \INT|\nabla u|^{p-2}\nabla u\nabla(u_n-u)dx\rightarrow0,\quad \INT h(x)|u|^{s-2}u(u_n-u)dx\rightarrow0.
\end{eqnarray}
By Lemma \ref{lem2.1}, it follows that  $u_n\rightarrow u$ in $L^r_k(\RT)$ and $L^p_{g^{\pm}}(\RT)$, where $g^{\pm}(x)=\max\{\pm g(x),0\}$.
Thus, by H\"{o}lder inequality, we have
\begin{eqnarray}\label{44}
\left.
  \begin{array}{ll}
\displaystyle\INT k(x)|u|^{r-2}u(u_n-u)dx\rightarrow 0,\\[1em]
\displaystyle\INT g(x)|u|^{p-2}u(u_n-u)dx\rightarrow 0,
  \end{array}
\right.
\end{eqnarray}
and
\begin{eqnarray}\label{45}
\left.
  \begin{array}{ll}
\displaystyle\INT k(x)(|u_n|^{r-2}u_n-|u|^{r-2}u)(u_n-u)dx\rightarrow 0,\\[1em]
\displaystyle\INT\ g(x)(|u_n|^{p-2}u_n-|u|^{p-2}u)(u_n-u)dx\rightarrow 0.
  \end{array}
\right.
\end{eqnarray}
Using the Sobolev embedding $D^{1,q}(\RT)\rightarrow L^{q^*}(\RT)$ and H\"{o}lder inequality, we obtain
\begin{eqnarray*}
\begin{aligned}
    &\left| \INT a(x)\phi_u|u|^{m-2}u(u_n-u)dx \right|\\
   \le&\left(\INT (a(x))^{\frac{q^*}{q^*-1}}|u|^{(m-1)\frac{q^*}{q^*-1}}|u_n-u|^{\frac{q^*}{q^*-1}}dx\right)^{\frac{q^*-1}{q^*}}\left( \INT |\phi_u|^{q^*}dx  \right)^{\frac{1}{q^*}}.
  \end{aligned}
\end{eqnarray*}
Applying  H\"{o}lder inequality again, there exists $t_0>1$  such that,
\begin{eqnarray*}
\begin{aligned}
    &\left|\INT a(x)^{\frac{q^*}{q^*-1}}|u|^{(m-1)\frac{q^*}{q^*-1}}|u_n-u|^{\frac{q^*}{q^*-1}}dx\right|\\
    \le& \left( \INT a(x)|u_n-u|^{t_0{\frac{q^*}{q^*-1}}}dx \right)^{\frac{1}{t_0}}
   \left(\INT (a(x))^{(\frac{q^*}{q^*-1}-\frac{1}{t_0})\frac{t_0}{t_0-1}}|u|^{(m-1)\frac{q^*}{q^*-1}\frac{t_0}{t_0-1}}dx  \right)^{\frac{t_0-1}{t_0}},
   \end{aligned}
\end{eqnarray*}
where
\begin{center}
   $1<t_0<\frac{p^*(q^*-1)}{q^*}$ and $1<\frac{q^*}{q^*-1}(m-1)\left(1+\frac{1}{t_0-1}\right)<p^*.$
\end{center}
Define the function
$$f(t)=\frac{q^*}{q^*-1}(m-1)\left(1+\frac{1}{t-1}\right).$$
Clearly, $f$ is decreasing in $t$. For $1 < t < \frac{p^*(q^* - 1)}{q^*}$,
$$f(1)=+\infty.$$
Since  $1<m<\frac{p^*(q^*-1)}{q^*}$, then
\begin{eqnarray*}
    0<f(\frac{p^*(q^*-1)}{q^*})=\frac{(m-1)p^*q^*}{p^*(q^*-1)-q^*}<p^*.
\end{eqnarray*}
By the intermediate value theorem, there exists a $t_0\in (1,\frac{p^*(q^*-1)}{q^*})$  such that $1<f(t_0)<p^*$.
Hence,
\begin{eqnarray}\label{46}
    \INT a(x)\phi_u|u|^{m-2}u(u_n-u)dx\rightarrow 0.
\end{eqnarray}
Moreover,
\begin{eqnarray}\label{48}
\INT a(x)\left(\phi_{u_n}|u_n|^{m-2}u_n-\phi_{u}|u|^{m-2}u\right)(u_n-u)dx\rightarrow 0.
\end{eqnarray}
From \eqref{43}, \eqref{44} and \eqref{46}, we have
$$\langle I'(u),u_n-u\rangle\rightarrow 0.$$
Note that $I'(u_n)\rightarrow 0$, there holds
$$\langle I'(u_n)-I'(u),u_n-u\rangle\rightarrow 0.$$
On the other hand,
\begin{eqnarray}\label{47}
\begin{aligned}
&\langle I'(u_n)-I'(u),u_n-u\rangle \\
\ge&\left[\left(\INT |\nabla u_n|^pdx\right)^{\frac{p-1}{p}}-\left(\INT |\nabla u|^pdx\right)^{\frac{p-1}{p}}\right]\\
&\left[\left(\INT |\nabla u_n|^pdx\right)^{\frac{1}{p}}-\left(\INT |\nabla u|^pdx\right)^{\frac{1}{p}}\right]\\
&+\left[\left(\INT h(x)| u_n|^sdx\right)^{\frac{s-1}{s}}-\left(\INT h(x)| u|^sdx\right)^{\frac{s-1}{s}}\right]\\
&\left[\left(\INT h(x)| u_n|^sdx\right)^{\frac{1}{s}}-\left(\INT h(x) | u|^sdx\right)^{\frac{1}{s}}\right]\\
&-\lambda\INT g(x)\left(|u_n|^{p-2}u_n-|u|^{p-2}u \right)(u_n-u)dx\\
&-\INT k(x)\left(|u_n|^{r-2}u_n-|u|^{r-2}u \right)(u_n-u)dx\\
&+\INT a(x)\left(\phi_{u_n}|u_n|^{m-2}-\phi_{u}|u|^{m-2} \right)(u_n-u)dx.
\end{aligned}
\end{eqnarray}
For a real function $f_b:\mathbb{R}^+\rightarrow \mathbb{R}$ defined by  $f_b(t)=(t^{b-1}-\alpha^{b-1})(t-\alpha)$ with $b>1$, we always have $f_b(t)\ge 0$. Thus, from \eqref{45}, \eqref{48}, \eqref{47} and $\langle I'(u_n)-I'(u),u_n-u\rangle\rightarrow 0$, it follows that
\begin{eqnarray}\label{49}
\begin{aligned}
&\left[\left(\INT |\nabla u_n|^pdx\right)^{\frac{p-1}{p}}-\left(\INT |\nabla u|^pdx\right)^{\frac{p-1}{p}}\right]\\
&\times\left[\left(\INT |\nabla u_n|^pdx\right)^{\frac{1}{p}}-\left(\INT |\nabla u|^pdx\right)^{\frac{1}{p}}\right]\rightarrow 0,
\end{aligned}
\end{eqnarray}
\begin{eqnarray}\label{410}
\begin{aligned}
&\left[\left(\INT h(x)| u_n|^sdx\right)^{\frac{s-1}{s}}-\left(\INT h(x)| u|^sdx\right)^{\frac{s-1}{s}}\right]\\
&\times\left[\left(\INT h(x)| u_n|^sdx\right)^{\frac{1}{s}}-\left(\INT h(x) | u|^sdx\right)^{\frac{1}{s}}\right] \rightarrow 0 .
\end{aligned}
\end{eqnarray}
Note that if $f_b(t_n) \to 0$, then $t_n \to \alpha$.
Thus, from \eqref{49} and \eqref{410},  we obtain
\begin{eqnarray*}
\left(\INT|\nabla u_n|^pdx\right)^{\frac{1}{p}}\rightarrow\left(\INT|\nabla u|^pdx\right)^{\frac{1}{p}},\\  \left(\INT h(x)|u_n|^sdx\right)^{\frac{1}{s}}\rightarrow\left(\INT h(x)|u|^sdx\right)^{\frac{1}{s}}.
\end{eqnarray*}
That is, $\|u_n\|_{D^{1,p}(\RT)}\rightarrow\|u\|_{D^{1,p}(\RT)}$ and $\|h^{\frac{1}{s}}u_n\|_{s}\rightarrow\|h^{\frac{1}{s}}u\|_{s}$. Combining with
\begin{center}
  $u_n\rightharpoonup u$ in $D^{1,p}(\RT)$,\quad $(h(x))^{\frac{1}{s}}u_n\rightharpoonup (h(x))^{\frac{1}{s}}u$ in $L^s(\RT)$,
\end{center}
the uniformly convexity of $D^{1,p}(\RT)$ and $L^s(\RT)$ implies $u_n\rightarrow u$ in $D^{1,p}(\RT)$ and $(h(x))^{\frac{1}{s}}u_n\rightarrow (h(x))^{\frac{1}{s}}u$ in $L^s(\RT)$. Therefore,
\begin{eqnarray*}
\begin{aligned}
 \|u_n-u\|=&\|u_n-u\|_{D^{1,p}(\RT)}+|u_n-u|_{s,h}\\
 =&\|u_n-u\|_{D^{1,p}(\RT)}+\|h^{\frac{1}{s}}u_n-h^{\frac{1}{s}}u\|_{s}\rightarrow0.
 \end{aligned}
\end{eqnarray*}
Thus, $u_n\rightarrow u$ in $E$. We complete the proof.
\end{proof}

Having verified the (PS) condition, we now investigate the geometry of $I$. First, note that $I$ is clearly even.
Let
$$\Omega:=\{ x\in\mathbb{R}^3|k(x)=0\}$$
and
$$X:=\{u\in E|u(x)=0 \text { a.e. } x \in \Omega\}.$$
Then $X$ is an infinite-dimensional linear subspace of $E$. We introduce a seminorm $\left[\cdot\right]_r$ on $X$ by
$$[u]_r:=\left(\INT k(x)|u|^rdx \right)^{\frac{1}{r}}.$$
\begin{lem}
  The seminorm $\left [ \cdot  \right ] _r$ is a norm on $X$.
\end{lem}
\begin{proof}
It suffices to show that
$$u\in X,\quad [u]_r=0\Rightarrow u=0.$$
Since $k(x)\ge0$, we have
$$0=[u]^r_r=\INT k(x)|u|^rdx=\int_{\{x\in\mathbb{R}^3| k(x)>0\}} k(x)|u|^rdx.$$
This implies $u=0$ a.e. on $\{ x\in\mathbb{R}^3|k(x)>0\}$. But $u\in X$, that is $u=0$ a.e. on $\{ x\in\mathbb{R}^3|k(x)=0\}$.
 Therefore,  $u=0$ a.e. on $\mathbb{R}^3$.
\end{proof}

Let $\Sigma$ denote the class of closed symmetric subsets of $E \setminus\{0\}$. For $A\in \Sigma$,  the genus $\gamma(A)$ is denoted by
$$\gamma(A):=\min\{\theta\in\mathbb{N}| \exists\varphi\in C(A,\mathbb{R}^3\setminus\{0\} ), \text{such that } \varphi(x)=-\varphi(-x)\}.$$
If no such minimum exists, we set $\gamma(A)=+\infty$. For properties of the genus, we refer to \cite{1986Chang,1986Rabinowitz,1990Struwe,1996Willem},  we omit them here.
\begin{lem}\label{lem4.4}
Let $\theta \in\mathbb{N}$, there is $\varepsilon=\varepsilon(\theta)$ such that
$$\gamma(\{u\in E|I(u)\le-\varepsilon\})\ge \theta.$$
\end{lem}
\begin{proof}
It follows from H\"{o}lder inequality that
$$-\frac{\lambda}{p}\INT g(x)|u|^pdx\le \left|\frac{\lambda}{p}\INT g(x)|u|^pdx \right|\le\frac{|\lambda|}{p}\|g\|_{\frac{3}{p}}\|u\|^{p}_{p^*}. $$
Since $\dim X=\infty$, for $\theta\in\mathbb{N}$, let $X_\theta$ be a $\theta$-dimensional subspace of $X$. Then for $u\in X_\theta$, by Lemma \ref{lem2.2} (i) and Proposition \ref{Pro2} (iii), we obtain
\begin{eqnarray}
\begin{aligned}
I(u)&\le C(\|u\|^p+\|u\|^s)-\frac{\lambda }{p}\INT g(x)\left | u \right |^pdx-\frac{1 }{r} \INT k(x)  \left |u  \right |^{r}dx +\frac{q-1}{qm}\INT a(x)\phi_{u} \left |  u\right | ^{m} dx\\
&\le C\|u\|^p+C\|u\|^s+\frac{|\lambda|}{p}\|g\|_{\frac{3}{p}}\|u\|^{p}_{p^*}-\frac{1}{r}[u]^{r}_r+C\|u\|^{\frac{qm}{q-1}}.
\end{aligned}\label{n4.11}
\end{eqnarray}
Since  $[\cdot]_r$ is a norm on $X$, and all norms on the finite-dimensional space $X_\theta$ are equivalent.  For $u\in X_\theta$, by \eqref{n4.11}, we deduce that
\begin{eqnarray}
I(u)\le C\|u\|^p+C\|u\|^s+B\|u\|^p-D\|u\|^r+K\|u\|^{\frac{qm}{q-1}},\label{n4.12}
\end{eqnarray}
where $ B, D, K>0$ are constants.  Note that $r<p<s$ and $r<\frac{qm}{q-1}$, from \eqref{n4.12}, we can choose $\varepsilon>0$(which depends on $\theta$) and $\rho>0$ such that $I(u)\le-\varepsilon$, if $\|u\|=\rho$.

Let $S:=\{u\in X_\theta| \|u\|\ = \rho\}$, then we have
$$S\subset\{u\in E|I(u)\le-\varepsilon\},$$
by the properties of genus, we get that
$$\gamma(\{u\in E|I(u)\le-\varepsilon\})\ge\gamma(S)=\theta.$$We complete the proof.
\end{proof}

Let
$$\Sigma_\theta:=\{A\in \Sigma|\gamma(A)\ge \theta \}.$$
From Lemma \ref{lem4.4}, we define a sequence of minimax values
$$c_{\theta}:=\inf _{A \in \Sigma_{\theta}} \sup _{u \in A} I(u).$$
Clearly,
$$c_1\le c_2\le \cdot\cdot\cdot\le c_\theta\le c_{\theta+1}\le \cdot\cdot\cdot .$$
Furthermore, by Lemma \ref{lem4.1},\ $I$ is coercive, thus bounded from below. That is, $c_\theta>-\infty$ for $\forall \theta\in \mathbb{N}$.

 For $c\in \mathbb{R}$, we define $K_c:=\{u\in E|I(u)=c,I'(u)=0\}$. By a standard argument (see \cite{1991AzoreroAlonso}), we have the following result.

\begin{lem}\label{lem4.5}
 Each $c_\theta$ is critical value  of $I$. Moreover, if
 $$c=c_\theta=c_{\theta+1}=\cdot \cdot\cdot=c_{\theta+\tau},$$
  then $\gamma(K_c)\ge \tau+1$.
\end{lem}

\begin{proof}[Proof of \textup{Theorem \ref{thm1.1}}]
For $\forall \theta\in \mathbb{N} $,   Lemma \ref{lem4.4}  implies the existence of $\varepsilon(\theta)>0$  such that
\begin{eqnarray*}
\{u\in E|I(u)\le-\varepsilon(\theta)\}\subset\Sigma_\theta,
\end{eqnarray*}
 thus $c_\theta\le -\varepsilon(\theta)<0$. By Lemma \ref{lem4.5}, there are infinitely many critical points of   $I$ with negative critical values. Thus,  problem \eqref{equ} admits infinitely many solutions with negative energy.  This completes the proof of Theorem \ref{thm1.1}.
\end{proof}

\section*{Acknowledgment} \     This work is supported by NSFC12301144 and  Sichuan Science and Technology Program (2024NSFSC1342).

\end{document}